\documentclass{amsart}

\usepackage{amssymb}
\usepackage{amsfonts}
\usepackage[all]{xy}
\usepackage[dvips]{graphics}

\def\mystretch{1.2}
\def\baselinestretch{\mystretch}

\def\baselinestretch{\mystretch}

\newtheorem{theorem}{Theorem}

\newtheorem{proposition}[theorem]{Proposition}
\newtheorem{lemma}[theorem]{Lemma}

\newtheorem{dfn}[theorem]{Definition}

\newtheorem{remark}[theorem]{Remark}

\theoremstyle{remark}

\def\abstract{\noindent \vspace{0.3cm}\large{\bf Abstract}\\
\small\def\baselinestretch{1}\normalsize}

\newcommand{\field}{{{\mathbb K}}}
\newcommand{\Alt}{{\mbox{Alt}}}
\newcommand{\Hom}{{\mbox{Hom}}}
\newcommand{\had}{{\mbox{\textbf{iad}}}\ }

\newcommand{\La}{{\mathcal {L}}}
\newcommand{\QB}{{\mathcal {QB}}}
\newcommand{\A}{{\mathbb {A}}^\hbar}
\newcommand{\C}{{\mathbb {C}}}

\newcommand{\Z}{{\mathbb {Z}}}
\newcommand{\N}{{\mathbb {N}}}
\newcommand{\R}{{\mathbb {R}}}
\newcommand{\Sym}{{Sym}}

\newcommand{\g}{{\mathfrak g}}

\newcommand{\V}{{\widetilde{V}}}
\newcommand{\Vast}{{{V}^\ast}}
\newcommand{\Vstar}{{{\widetilde{V}^\ast}}}
\newcommand{\ad}{{\tt{ad}}}

\newcommand{\B}{{\mathcal B}}

\begin{document}

\title{Strongly homotopy Lie  bialgebras and Lie quasi-bialgebras}

\author {Olga Kravchenko }
\address{Universit\'e de Lyon,
Universit\'e Lyon1, CNRS, UMR 5208 Institut Camille Jordan, Batiment
du Doyen Jean Braconnier, 43, blvd du 11 novembre 1918, F - 69622
Villeurbanne Cedex, France} \email{ok@alum.mit.edu}

\date{April 27, 2007}

\dedicatory{This paper is dedicated to Jean-Louis Loday on the
occasion of his 60th birthday with admiration and gratitude}

 \vspace{1cm}

\maketitle

\thispagestyle{empty}

\noindent\begin{abstract} \noindent Structures of Lie algebras, Lie
coalgebras, Lie bialgebras and Lie quasibialgebras are presented as
solutions of Maurer-Cartan equations on corresponding governing
differential graded Lie algebras using the big bracket construction
of Kosmann-Schwarzbach. This approach provides a definition of an
$L_\infty$-(quasi)bialgebra (strongly homotopy Lie
(quasi)bialgebra). We recover an $L_\infty$-algebra structure as a
particular case of our construction. The formal geometry
interpretation leads to a definition of an $L_\infty$
(quasi)bialgebra structure on $V$ as a differential operator $Q$ on
$V,$ self-commuting with respect to the big bracket.   Finally, we
establish an $L_\infty$-version of a Manin (quasi) triple  and get a
correspondence theorem with $L_\infty$-(quasi) bialgebras.
\end{abstract}

\noindent

\section{Introduction.}

 Algebraic structures are often  defined  as certain
maps which must satisfy quadratic relations. One of the examples is
a Lie algebra structure: a Lie bracket satisfies the Jacobi identity
(indeed the Jacobi identity is a quadratic relation since the
bracket appears twice in each summand).  Other examples include an
associative multiplication (the associativity condition is
quadratic), $L_\infty$ and $A_\infty$ algebras (also called strongly
homotopy Lie and strongly homotopy  associative algebras) and many
others.

The subject of this article is a description of Lie (quasi)
bialgebras and their $L_\infty$-versions. The main philosophy is
that the axioms of Lie (quasi)bialgebras (and their
$L_\infty$-versions) could be written in the form of a quadratic
relation on a certain governing differential graded Lie algebra. We
find the governing differential graded Lie algebras for Lie
bialgebra and Lie quasi-bialgebra structures using
Kosmann-Schwarzbach's big-bracket construction \cite{KS}. The
$L_\infty$ brackets are obtained by using the (higher) derived
brackets  \cite{KS4,Vo,AI}.

The quadratic relation on the structure can be expressed in the form
of a Maurer-Cartan equation.

Classically, solutions of the Maurer-Cartan equation are considered
only from the first graded component of the governing
 differential graded Lie algebra and they give the original algebraic
 structure.
  However solutions of the Maurer-Cartan equation  in the whole governing
 differential graded Lie algebra provide a strongly homotopy version of the
original one.

Let us give here the definitions of Lie bialgebras, Lie
quasibialgebras and Manin triples and pairs to start with.

A good reference on  Lie (quasi)bialgebras is the book by Etingof
and Schiffmann (\cite{ESch}, pages 32--34 and 150--152).

\begin{dfn} \label{biLie} A Lie bialgebra structure on a vector space
$V$ is the following data:
\begin{description}
  \item[a] a Lie bracket, $\{ \cdot, \cdot \} : {V} \wedge {V} \to  V;$
  \item[b] a Lie cobracket, that is an
element $\delta: {V} \to {V} \wedge {V},$ satisfying the coJacobi
identity: $$\Alt(\delta \otimes 1)\delta(x) = 0$$
  \item[c]  a compatibility condition between
  $ \{\cdot, \cdot \}$ and $\delta,$ meaning that  $\delta$ is a
1-cocycle: $\delta(\{x,y \}) = \{\delta(x), 1 \otimes y + y \otimes
1\} + \{1 \otimes x + x \otimes 1, \delta(y)\},$
\end{description}
\end{dfn}

 where  for $x,y,z$ elements of $V$ we denote $ Alt(x
\otimes y \otimes z) = x \otimes y \otimes z + y \otimes z \otimes x
+ z \otimes x \otimes y.$

To distinguish a Lie bracket on the governing differential graded
Lie algebra from a Lie bracket in a Lie (quasi)bialgebra structure
we denote the former one by $[\cdot, \cdot ]$ and the latter one by
$\{ \cdot, \cdot \} $ throughout the article.

\begin{dfn} \label{qbiLie} A Lie quasibialgebra structure on a vector space
$V$ is the following data:
\begin{description}
  \item[a] a Lie bracket $\{ \cdot, \cdot \};$
  \item[b] an
element $\delta \in \Hom({V}, {V} \wedge {V}),$ and an element $\phi
\in V \wedge V \wedge V$ satisfying a modified coJacobi identity:
$$\frac{1}{2}\Alt(\delta \otimes 1)\delta(x) = \{x \otimes 1 \otimes
1 + 1 \otimes x \otimes 1 + 1 \otimes 1 \otimes x, \phi\}$$ and
$(\delta \otimes 1\otimes 1 + 1 \otimes \delta \otimes 1 + 1 \otimes
1\otimes \delta )(\phi) = 0.$
  \item[c]   $\delta$ is a 1-cocycle with respect to the bracket $\{ \cdot, \cdot \}$.
\end{description}
\end{dfn}

The main difficulty in  defining strongly homotopy versions of Lie
bialgebras and Lie quasibialgebras in the generalized Maurer-Cartan
approach is in finding the corresponding governing Lie algebra. Let
$\bigwedge V$ be the vector space be the exterior power vector
space. Then the Lie bracket and the Lie cobracket on $V$ could be
considered as elements of the space of homomorphisms  $Hom(\bigwedge
V, \bigwedge V).$ However, the Lie algebra structure on this space
$Hom(\bigwedge V, \bigwedge V)$ is not obvious. Actually, its
associative analogue, $Hom(\otimes V, \otimes V),$ does not have
such a structure at all (which makes the Lie bialgebra quantization
problem so intriguing).

We define the correct Lie algebra structure on $Hom(\bigwedge V,
\bigwedge V)$ by considering an isomorphism of this space to
$\bigwedge V^\ast \otimes  \bigwedge V.$ The Lie algebra structure
on the latter  is  given by the Kosmann-Schwarzbach big bracket
(\cite{KS}), which is also just the  Batalin-Vilkovisky  odd Poisson
bracket on the space $V^\ast \otimes V$ with a shifted degree
(\cite{BV,FV}).

The isomorphism from  $\bigwedge V^\ast \otimes  \bigwedge V$ to
$Hom(\bigwedge V, \bigwedge V)$ is obtained by using the iterated
adjoint Lie action. In the part``$L_\infty$-structures'' we first
interpret the strongly homotopy Lie algebra structure in terms of
the iterated adjoint action and then generalize the construction  to
the case of Lie (quasi)bialgebras. This might look as a complicated
way to define an
 $L_\infty$- structure while
 a definition by coderivations of $\bigwedge\nolimits^k V$ is more conventional.
However, a Lie (quasi)-bialgebra structure  is not given by a
coderivation. The higher derived brackets introduced by
 Th.~Voronov (\cite{Vo}) and also studied by Akman-Ionescu
(\cite{AI}) provide  the necessary tool.

Finite dimensional Lie bialgebras are in one-to-one correspondence
with Manin triples:

\begin{dfn} A Manin triple $(\g, \g_+, \g_-)$ is a triple of finite
dimensional Lie algebras, where $\g_+  \bigoplus \g_- = \g$ as a
vector space and $\g$ is equipped with a nondegenerate symmetric
invariant bilinear form $< \cdot, \cdot>$ such that   $\g_+ $ and
$\g_-$ are Lagrangian subalgebras (that is maximal isotropic
subspaces which are  Lie subalgebras).
\end{dfn}

 Lie quasi-bialgebras turn out to be described by the notion of
a Manin pair.
\begin{dfn} A Manin pair is a pair $(\g, \g_+)$ where $\g$ is a finite
dimensional Lie algebra equipped with a non-degenerate symmetric
invariant bilinear form $< \cdot, \cdot>$ and   $\g_+ $ is a
Lagrangian subalgebra.

A  Manin quasi-triple (also called a marked pair) is a pair $(\g,
\g_+)$ with a chosen Lagrangian complement of $\g_+.$
\end{dfn}

The main theorem here is by Drinfeld \cite{D} which states that
 Manin quasi-triples are in one-to-one correspondence with  Lie
quasi-bialgebras.

 We will formulate and prove an $L_\infty$ (in other words strongly homotopy) version
 of the Lie (quasi)bialgebras - Manin (quasi)triples correspondence.

It should be mentioned that an operad \cite{MSS} or rather properad
(or PROP) approach is not used in this paper. However, the
definition of an $L_\infty$-bialgebra coming from a minimal
resolution of a Lie bialgebra PROP coincides with ours as shown
explicitly in the work of Sergei Merkulov \cite[Corollary 5.2]{M1}
and also could be derived from works \cite{Gan} on dioperads and
\cite{MV, Val} on Koszul PROPs, the Lie bialgebra PROP being one of
them.

To avoid confusion, it should also be mentioned that if one wants to
consider an associative bialgebra there is no similar construction
of a governing Lie algebra since the bialgebra PROP is not Koszul.
Our method works only for Lie (quasi)bialgebras.

We use the Koszul sign convention: in a graded algebra whenever
there is a change of places of two symbols there should be a
corresponding sign. Throughout this paper, the summation convention
is understood: indices $\alpha,\beta,\dots$ once as superscript and
once as subscript in a formula are to be summed over.

\subsection*{Acknowledgements}
My interest in homotopy structures flourished when I was a postdoc
in Strasbourg with Jean-Louis Loday. One of his inventions,  his
Leibniz algebras, I encounter over and over again, and in particular
in this article. I consider myself extremely lucky to have profited
from his amazing intellectual generosity and his joy in mathematics.

 I thank Ezra Getzler for useful discussions. I am grateful to
Vladimir Hinich who helped me to see cohomological implications and
encouraged me to write them up. I am indebted to Boris Tsygan for a
discussion  crucial to the completion  of this article (in
particular for Remark \ref{thm}) and also for making my work look
more exciting by pointing out possible connections to \cite{EGH}. I
thank Sergei Merkulov for pointing out references \cite{MV, M1} and
answering my PROPeradic questions.  The article  also benefited from
comments of Jim Stasheff. Last but not least, I cannot overestimate
the help of Yvette Kosmann-Schwarzbach in putting  this article in
its final shape.

I am grateful to the referee for thorough reading and useful
remarks.

\section{Kosmann-Schwarzbach's big bracket }

Treating  exterior powers of a sum of a vector space with its dual
as a super-Poisson algebra was pioneered in the works of
Batalin-Fradkin-Vilkovisky \cite{FV,BV} and later of Stasheff
\cite{S}  and Kostant-Sternberg \cite{KoSt}. In 1991 Yvette
Kosmann-Schwarzbach published an article \cite{KS} where the term
\emph{big bracket} was introduced in order to describe
proto-bialgebras (a notion generalizing Lie bialgebras). This big
bracket defines, in particular, the Lie structure of the governing
dgLie algebra of Lie bialgebras (and, in fact, Lie algebras, Lie
coalgebras and Lie quasibialgebras).

Here is the construction from \cite{KS} in a $\Z$-graded context. On
a direct sum of a vector space and its dual there is a symmetric
bracket pairing the space with its dual. There is a way to see this
symmetric bracket as a graded odd Lie bracket on a graded space.

 On any $\Z$-graded space $X = \sum X^{i}$ there is an operation
called de-suspension, mapping $X$ to the space $X[1]$ so that
$X[1]_i = X_{i+1}. $ This shift of the degree allows in particular
to see symmetric powers of the  space $X$ as exterior powers of the
shifted space $X[1]$ (and vice-versa). There is the following
natural identification of symmetric and exterior powers of $X:$
\begin{equation} \label{sym}
\Sym^n (X) \simeq \left(\bigwedge\nolimits^n (X[1])\right)[-n].
\end{equation}
Thus a symmetric form acting from $\Sym^2(X)$ to $\field =
\Sym^0(X)$ becomes a  bracket
$$ \left(\bigwedge\nolimits^2(X[1])\right)[-2]
\ \to \ \bigwedge\nolimits^0(X[1])= \field.$$  This bracket could
also be viewed as a map: $\bigwedge\nolimits^2(X[1])\to
 \left(\bigwedge\nolimits^0(X[1])\right)[2]= \field[2].$ If $X$ is ungraded we get
elements of degree $-1$ in $X[1]$ and the $\field[2]$ is in degree
$-2.$ Thus the bracket has degree $0,$ since two elements of degree
$-1$ are send by the bracket to $\field[2]$ of degree $(-1) + (-1) =
-2.$

Let $\widetilde{V}$ be a finite dimensional $\field$-vector space
($\field=\R$ or $\C$), and $\widetilde{V}^\ast$ its dual.
 A non-graded space
$\widetilde{V}^\ast \oplus \widetilde{V}$ could be considered as a
graded one by assigning degree $0$ to each of its elements. Consider
$\Vast \oplus V  = (\widetilde{V}^* \oplus \widetilde{V})[1],$ this
means that ``points'' of $\Vast \oplus V$ are in degree $-1.$
Algebraic  functions on $\V \oplus \Vstar$ form the space
\begin{equation}\label{B}
B =  {\oplus_{j \geq -2} B^{j},} \quad \text{where} \ B^{j} =
{\oplus_{p+q = j} \left(\bigwedge\nolimits^{p+1}\Vast \otimes
\bigwedge\nolimits^{q+1}V\right), \ j \geq -1}, \quad \text{and} \
B^{-2} = \field.
\end{equation}
  From this point of view $B$ is an
algebra of exterior powers  of the odd space $\Vast \oplus V.$

Grading in $B$ is  given by the sum $p+q.$ That is, $B^j$ consists
of terms with $j = p+q,$ for $ p,q \geq -1,$ so that the first few
terms are as follows:
\begin{eqnarray*}
B^{-2} & = &k,\\
B^{-1} &= & V \quad \oplus \quad \Vast,\  \\
 B^{0} \ \ & = & V \wedge V    \quad \oplus  \quad \Vast
\otimes V \quad \oplus \quad \Vast \wedge \Vast,\\
B^1 \ \ & =  &V \wedge V \wedge V \quad \oplus \quad \Vast \otimes V
\wedge V
 \quad  \oplus  \quad V^\ast \wedge V^\ast \otimes V \quad \oplus \quad V^\ast
 \wedge V^\ast \wedge V^\ast.
 \end{eqnarray*}
The space $B$ is bigraded as follows: \vspace{5mm}

\begin{tabular}{c|ccccc}

  %  after  \\:  \hline or \cline{col1-col2} \cline{col3-col4} ...
& $\cdots$  &     &     &    & \\
   & &     &     &    & \\
 2 \qquad & $V^\ast    \wedge     V^\ast    \wedge     \Vast$
  & $\cdots$     &     &  &   \\
   & &     &     &    & \\
 1 &$V^\ast  \wedge   \Vast$   &  $V^\ast \wedge  V^\ast \ \otimes  \ V $
 &  $\cdots$ &  &\\
  &&     &     &   &  \\
  0 & $V^\ast $  & $V^\ast\ \otimes \ V $ & $V^\ast\ \otimes \ V  \wedge V $
  &$\cdots$  &\\
   &&     &     &   &  \\
 -1 & $\field$ & $ V $ & $V  \wedge V $ & $V  \wedge V  \wedge
   V $ & $\cdots$ \\
   & &     &     &    & \\ \hline
 \ $p \quad \diagup \quad q$&\qquad -1 &\qquad 0   &  \qquad 1 &  \qquad 2 & \\

\end{tabular}

\vspace{5mm}

Let $\langle \cdot, \cdot \rangle$ be the natural pairing of $\V$
and $\widetilde{V}^*.$ We extend it to a symmetric form on $\V
\oplus \Vstar$ as follows: for $x,y  \in \V $ and $v,w \in \Vstar:$
$$ \langle x+v, y+w \rangle =\langle x, w \rangle + \langle v, y
\rangle.$$ This symmetric form on $\V \oplus \Vstar$ could be
considered as an antisymmetric form on the de-suspended space $V
\oplus \Vast.$ Moreover, on $B$ it gives a Lie algebra structure by
the following

\begin{dfn}
\emph{The big bracket} is the graded Lie algebra structure on
algebraic functions on $V \oplus V^\ast$ defined as follows.
\begin{itemize}
\item For $ u,v \in B^{-2} \oplus B^{-1} = k \oplus V \oplus
V^\ast:$
$$[u,v] = \left\{
\begin{array}{ll}
\langle u, v \rangle & \mbox{if} \quad  u \in V \oplus
V^\ast \ \mbox{and} \ v \in V \oplus V^\ast\\
0& \mbox{if} \quad  u \in k  \quad \mbox{or} \quad  v \in k
\end{array}
\right.
$$
\item The bracket  on other terms is defined by linearity and the
graded Leibniz rule: for $u \in B^{k},  \ v  \in B^{l}, \ w \in
B^{m}$
$$
[u,v \wedge w] = [u,v] \wedge w + (-1) ^{kl}v \wedge [u,w]
$$
\end{itemize}
\end{dfn}
\begin{remark}\label{sub}

\begin{enumerate}
\item $[u,v] = - (-1)^{kl}[v,u]; \quad u \in B^{k},  \ v  \in
 B^{l},$ that is the big bracket is skew-symmetric in the graded
 sense;
 \item  $[\cdot, \cdot]: B^i \wedge B^j \to
B^{i+j},$  in other words: the bracket is of degree zero;
\item in particular,
$$[\cdot, \cdot]: \ (\bigwedge\nolimits^k \Vast \otimes \bigwedge\nolimits^l V) \
\wedge \  (\bigwedge\nolimits^m \Vast \otimes \bigwedge\nolimits^n
V) \longrightarrow \bigwedge\nolimits^{k+m-1} \Vast \otimes
\bigwedge\nolimits^{l+n-1} V;$$
  \item $B^{0}$ is a Lie subalgebra of $B;$
  \item $B^{-2}$ is the center of $(B, [ \cdot, \cdot]).$
\end{enumerate}
\end{remark}

\section{Governing graded Lie algebras and Maurer-Cartan equations}

The graded Lie algebra $B$ has several graded Lie subalgebras.
Verification that they are indeed Lie subalgebras is easy using
Remark \ref{sub} (3). Some of these Lie subalgebras give well known
structures as  was described in \cite{KS}. We put them in the
following table containing the $-1, 0$ and $1$ graded components of
the governing Lie algebras: \vspace{5mm}

\noindent \begin{tabular}{|c||c|l|l|c|} \hline &&&& \\
$\g $ - Lie &&&& Solutions of\\
  subalgebra  &  $ \g^{-1} $& $\qquad \qquad \g^{0} $   & $\qquad \qquad \qquad\g^1 $
 &  Maurer-Cartan equation \\  of $B$ &&&& \\
  \hline \hline
 &&&& \\
  $\mathcal L:$ & $u \in V$& $f \in  V^\ast\ \otimes \ V$
  & $l \in V^\ast \wedge \Vast \ \otimes \ V$
  & $[l , l ] = 0$\\
 column $q =0$&&&& Lie algebra\\
  &&&& structure  on $V$\\
  \hline &&&& \\
 $\mathcal C:$  & $u^\ast \in \Vast$& $ f \in V^\ast\ \otimes \ V$
 & $ c \in V^\ast\ \otimes \ V \wedge V$
  & $[c,c] = 0$ \\
 row $ p=0$&&&&Lie coalgebra \\
 &&&& structure on $V$ \\
 \hline &&&& \\
 $\mathcal B:$    & --- &
 $ f \in V^\ast\ \otimes \ V $  & $ c \in
 \quad  V^\ast\ \otimes \ V  \wedge V \  $ &
    $[c+l , c +l ] = 0$ \\
    $p,q \not= -1$ &&&$
l \in \quad  V^\ast \wedge  V^\ast \ \otimes  \ V$
    & Lie bialgebra \\&&&& structure on $V$ \\
    \hline &&&& \\
 $\mathcal QB:$  & $u \in V$&
 $f \in \quad V^\ast\ \otimes \ V $  &
 $ c \in \quad V^\ast\ \otimes \ V  \wedge V
    \ $   & $[c + l + \phi,
    c+l+\phi] = 0$\\
 {$q \not= -1$ } & & $g \in \quad  \ V  \wedge V \ $ & $
l \in \quad \ V^\ast \wedge  V^\ast \ \otimes  \ V\ $ & Lie
    quasi-bialgebra  \\
   &&&$\
 \phi \in \quad V  \wedge V  \wedge
   V \ $& structure on $V$\\ \hline
\end{tabular}

\vspace{5mm}

Let $\bigwedge\nolimits V$ denote the sum $ k  \ \oplus V \oplus
\bigwedge\nolimits^2 V \oplus \bigwedge\nolimits^3 V  \oplus
\cdots.$
 Then
 $ {\mathcal C}  = \Vast \otimes \bigwedge\nolimits V \simeq Hom (V, \ \bigwedge\nolimits V),$ and
 ${\mathcal L} = \bigwedge\nolimits \Vast \otimes V \simeq Hom (\bigwedge\nolimits V, \ V).$
In other words,

\begin{proposition}\label{zoo}
On a vector space $V$
\begin{description}
 \item[Lie algebras ]
${\mathcal L}  = \oplus_{k\geq 0} \bigwedge\nolimits^k \Vast \otimes
V$ is \textit{the governing graded Lie algebra} of Lie algebra
structures on $V.$ In particular, an element of degree $1, l$ such
that $[l,l] = 0$ defines a Lie algebra structure $\{ \cdot , \cdot
\} \in Hom(V \wedge V, V)$ as follows: $ \{ x, y \} = [[l, x ], y
].$
\item[Lie coalgebras] ${\mathcal C} = \Vast \otimes (\oplus_{l\geq
0}\bigwedge\nolimits^l V) $ is \textit{the governing graded Lie
algebra} of Lie coalgebra structures on $V.$ In particular, an
element $c$ of degree $1$ such that $[c,c] = 0$ defines a Lie
coalgebra structure $\delta \in Hom(V, V \wedge V)$ as follows: $
\delta (x) = [c, x].$
  \item[Lie bialgebras]
${\mathcal B} = \oplus_{k\geq 1, l \geq 1} \bigwedge\nolimits^k
\Vast \otimes \bigwedge\nolimits^l V$ is \textit{the governing
graded Lie algebra} of Lie bialgebra structures on $V.$ In
particular, elements of degree $1$, $ \theta \in  V^\ast\ \otimes \
V  \wedge V $ and $l \in \ V^\ast \wedge V^\ast \ \otimes  \ V $
    such that $[ c +l , c+l  ]
= 0$ define a Lie bialgebra structure with the cobracket $ \delta
(x) = [c, x]$ and the bracket $ \{ x, y \} = [[l, x], y
].$
  \item[Lie quasi-bialgebras]
$\mathcal QB = \oplus_{k\geq 0, l \geq 0} \bigwedge\nolimits^k \Vast
\otimes \bigwedge\nolimits^l V$ is \textit{the governing graded Lie
algebra} of Lie quasibialgebra structures on $V.$ In particular,
elements of degree $1:$ $$\ c \in  V^\ast\ \otimes \ V  \wedge V, \
l \in\ V^\ast \wedge V^\ast \ \otimes  \ V  \ \mbox{and} \ \phi \in
V \wedge V \wedge V \quad  \mbox{such that}$$
\begin{equation} \label{quasi}
[c+l + \phi, c+ l + \phi] = 0
\end{equation}
    define a Lie quasibialgebra structure with the 1-cocycle $ \delta
(x) = [c, x],$ the Lie bracket $ \{ x, y \} = [[l, x ], y ]$ and the
3-tensor $\phi.$
\end{description}

\end{proposition}

\begin{proof}
Under the identification $\Hom(\bigwedge\nolimits^m V,
\bigwedge\nolimits^n V) \simeq \bigwedge\nolimits^m V^\ast \otimes
\bigwedge\nolimits^n V^\ast$ the condition $[l,l] = 0$ is exactly
the Jacobi identity.

The equation $[c,c]= 0$  is exactly the co-Jacobi identity on the
corresponding $\delta \in Hom (V, V \wedge V).$

The commutator $[l+c, l+c] $ lies in $ \left( \bigwedge\nolimits^3
\Vast \otimes V\right) \oplus \left( \Vast \otimes
\bigwedge\nolimits^3 V \right)\oplus \left( \bigwedge\nolimits^2
\Vast \otimes \bigwedge\nolimits^2 V \right)\subset B.$ Then $l+c $
defines a Lie bialgebra structure if its commutator with itself is
$0, $ that is
 all three components of  it in $ B^2$ must be equal to $0.$
 This leads to the three axioms of a Lie bialgebra:
\begin{description}
  \item[a] Jacobi identity follows from $[l,l] =0,$ for
  the derived bracket given by $  \{x, y\} = [[l, x] y].$
  \item[b] Co-Jacobi identity follows from $[c,c] = 0.$
  \item[c] The cocycle condition translates as: $[l,c] = 0.$
\end{description}
In fact, a presentation of a Lie bialgebra structure as a square
zero element in $\bigwedge\nolimits (V^\ast \oplus V )$ appeared
first in \cite{LR} (see also \cite{R} for the idea of a derived
bracket involved).

In the same manner the equation \ref{quasi} gives independent
equations:
 $[l,l] = 0, \quad [c,c] + 2[l, \phi]= 0, \quad
  \ [c, \phi] = 0 $ and
  $[c,l] = 0,$ which give  the axioms of a Lie quasibialgebra.

  One could look for complete proofs in \cite{KS}.
\end{proof}

 In all the cases of  Proposition \ref{zoo} we get a graded Lie
 algebra with a differential $d$ given by the adjoint action of
an auto-commuting element in the first degree. Hence $\mathcal C,
\mathcal L, \mathcal B  \ \mbox{and} \ \mathcal{QB}$ become dgLie
algebras. Here is a general fact:
\begin{proposition}{\cite{KS2}} \label{Leibniz}
A differential graded Lie algebra $(\g, \ [ \cdot, \cdot ], \ D)$
gives rise to a new bracket of degree $+1,$ called a derived bracket
with respect to $D:$ for $a,b \in \g, \ \{ a,b \} = [Da, b].$ This
is a Loday-Leibniz algebra bracket in the sense of  {\cite{L}}. If \
$W$ is an {\it{Abelian}} subalgebra of $\g,$ such that $[DW, W]$ is
in $W,$ then $\{\cdot, \cdot \}$ is a Lie bracket on $W$. Moreover,
 $D: \ (W, \{, \}) \to \ (\g, [,]) $ is a
Lie algebra morphism.
\end{proposition}

In our case, the algebra $\g$ is one of  $\mathcal C, \mathcal L,
\mathcal B  \ \mbox{and} \ \mathcal{QB},$ while the subalgebra Lie
$W$ in all cases is the same $B^{-1} = V^\ast \oplus V.$ The derived
bracket with respect to corresponding differentials on $\mathcal C,
\mathcal L, \mathcal B $ will define the Lie bracket on $V^\ast
\oplus V$, leading to  Manin triples
 (we could see a Lie structure on $V$ as a particular
case of a Manin triple with a zero cobracket, and a Lie coalgebra
structure as a Manin triple with zero  Lie bracket). In the same
way, we get a marked  Manin pair from $\mathcal{QB}.$

However, we need a certain refinement of  Proposition \ref{Leibniz},
since the Lie subalgebra $B^{-1}$ is not Abelian.

\begin{proposition} \label{hderived}
The differential graded Lie algebra $(B, \ [ \cdot, \cdot ], \ d),$
gives rise to a derived bracket of degree
 $1: \ \{ a,b \} = [da, b].$

This derived bracket restricted to $B^{-1}$ is a Lie bracket and $d:
\ (B^{-1}, \{\cdot, \cdot  \}) \to \ (B^{0}, [\cdot, \cdot ]) $ is a
Lie algebra morphism.
\end{proposition}

\begin{proof}
Notice that $\{ B^{-1}, B^{-1}\} \subset B^{-1}.$ The subspace
$B^{-1}= V^\ast\oplus V $ is not an Abelian Lie subalgebra of $(B,[
\cdot, \cdot ]),$ and we cannot use Proposition \ref{Leibniz}.
However, since the bracket $[ \cdot, \cdot ]$ on $B^{-1}$ takes
values in the center of $B$ (namely, $B^{-2} = k$), the Jacobi
identity on $B$ gives us the skew-symmetry of the derived bracket
$\{ \cdot, \cdot \}.$ Hence the derived bracket defines a Lie
algebra structure (not just Loday's) on $V^\ast\oplus V .$

The Lie algebra morphism part is immediate: $d\{a,b\} = d[da,b] =
[da,db].$
\end{proof}

Finally, to make explicit the connection to Manin pairs and triples
we have the following
\begin{proposition}
Let $l \in \bigwedge\nolimits^2 \Vast \otimes V, \ c \in \Vast
\otimes \bigwedge\nolimits^2 V, \ \phi \in \bigwedge\nolimits^3 V.$
Then $(V \oplus \Vast, \ \{ a,b \}= [da, b]),$  for
\begin{itemize}
   \item $ d = ad_{l+c},$ with a condition
   $[l +c, l+c] =0$
   defines a Manin triple $(V \oplus \Vast, V , \Vast).$
  \item  $ d = ad_{l + c +
\phi},$ satisfying $[{l + c + \phi},{l +c + \phi}] = 0$ defines  a
Manin pair  $(V \oplus \Vast, V).$
\end{itemize}
\end{proposition}

\section{$L_\infty$  structures}

In the previous section we have seen that a Lie algebra structure on
a space $V$ is obtained as a derived bracket defined on the graded
Lie algebra $\mathcal L$ with a differential given by an adjoint
action of an element  from  $V^\ast \wedge V^\ast
\otimes V \subset B^1$
 whose bracket with itself is $0.$

If we take an element from  $\La= \oplus_{k =0}^{\infty}
\bigwedge\nolimits^k V^\ast \otimes V$ (not just from $V^\ast \wedge
V^\ast \otimes V$) whose bracket with itself is $0$  we could define
an $L_\infty$
 structure on $V$ using  higher derived brackets (\cite{Vo},
\cite{AI}).  An autocommuting element from $\B$ (or $\mathcal{QB}$)
defines an $L_\infty$ bialgebra (or $L_\infty$ quasi-bialgebra)
structure using iterated adjoint action via derived  brackets. Here
we develop this theory.

A certain subtlety is in defining what an element of degree $1$
would mean in this context. For that we need to introduce a new
notion of degree.

\subsection{Degree}
 A starting point of any homotopy construction is a graded vector
space $(W = \oplus_a W^a,  \delta),$ with a differential of degree
$1,$   in other words  a complex
\begin{equation}
    \xymatrix{
\cdots \ar[r]^{\delta} & W^a  \ar[r]^{\delta} & W^{a+1}
\ar[r]^{\delta}& \cdots }.
  \end{equation}
  Let the space
$V$ be graded. Then the space $B = \oplus_{p,q \geq -1}
(\bigwedge\nolimits^{p+1}\Vast \otimes \bigwedge\nolimits^{q+1}V),$
as well as the space of  maps  $\bigwedge\nolimits V \to
\bigwedge\nolimits V,$ inherit the grading from $V$ in a consistent
way.
 We  define a grading and a differential on these spaces as follows.

  \begin{lemma} \label{dual-deg} Let  $(V = \oplus V^a, \delta)$
be a differential graded space over an ungraded  field $\field. $
%Assume that each $V_a$ is finite dimensional.
  Then on the dual space $\Vast$ there is a
corresponding grading with the opposite sign  and a differential $d$
of degree $1,   d: \Vast_{-(a+1)} \to \Vast_{-a}:$
\begin{equation}
    \xymatrix{
 \cdots & V_{-a}^\ast  \ar[l]_{d} &
V_{-(a+1)}^\ast \ar[l]_{d}& \ar[l]_{d} \cdots
\phantom{V_{-(a+1)}^\ast}}
  \end{equation}
  \end{lemma}
  \begin{proof}

 If $u$ is in $V^a$ and $v= \delta u \in V^{a+1}$ then
if $v^\ast \in \Vast$ is a dual of $v$ let us define: $u^\ast = d
v^\ast,$ where $d$ acts as an adjoint operator on the dual space.
This element  $d v^\ast $ is then dual to $u:$
\begin{equation}\label{dual d} 1=[v^\ast, v] = [v^\ast, \delta
u] = [d v^\ast, u] = [u^\ast,u].\end{equation}
 From the duality
of graded spaces $V$ and $V^\ast$ we get that $d$ is of degree $1$
on  $\Vast.$
 \end{proof}

 This way $(\Vast \oplus V, d+ \delta )$ becomes a
complex. We could extend the action of $d + \delta$ on $
\bigwedge\nolimits^p V^\ast \ \otimes \ \bigwedge\nolimits^q V$ for
any positive $p,q$ by the Leibniz rule, so that the whole $B$
becomes a complex.

 Let us define a new grading on the complex $B,$ taking into account
 the grading on $V:$
 \begin{dfn}\label{degree} Given a graded vector space $V = \oplus_{a}
 V^a$,
 consider elements from $ \bigwedge\nolimits^{p}
V^\ast \ \otimes \ \bigwedge\nolimits^{q} V.$ The degree of elements
from $V^\ast$ is defined according to Lemma \ref{dual-deg}, and
$\bigwedge\nolimits^0 V = \bigwedge\nolimits^0 V^\ast = \field$ is
of degree $0.$

Let $w = x^\ast_1 x^\ast_2...x^\ast_{p}y^1 y^2...y^{q} \in
\bigwedge\nolimits^{p} V^\ast \otimes \bigwedge\nolimits^{q} V, $
such that  $x^\ast_i, y^j$ are homogeneous, that is belonging to
just one graded component of $V^\ast \oplus V,$ the degree is
denoted by $\widetilde{\phantom{x_i}}.$ Then the \textit{internal
degree} of $w$ is the following sum of the degrees
 $\textbf{id}(w) = \sum \widetilde{x_i^\ast} + \sum
\widetilde{y^i}.$ The \textit{external degree} comes from the
grading of $B:
 \textbf{ed}(w) = p+q -2.$

 The \textit{total degree} of
$w$ then is $\textbf{td}(w) =  \textbf{id}(w) + \textbf{ed}(w) =
\sum \widetilde{x_i^\ast} + \sum \widetilde{y_i} + (p+q - 2).$

A map $\tau_{pq}: \bigwedge\nolimits^{p} V \to
\bigwedge\nolimits^{q} V$ sending $ v^1 v^2...v^{p}$ to $u^1
u^2...u^{q}$ has the degree $n$ if $(\sum \widetilde{u_i}) - (\sum
\widetilde{v_j}) + (p+q - 2) = n.$ \end{dfn}
\begin{remark} Let  $v^\ast_i$ be defined by the following property:
$[v^\ast_i,v^i] = 1.$ Then a map $\tau_{pq}$ has degree $n$ if and
only if  $t_{pq} = v^\ast_1 v^\ast_2...v^\ast_{p}u^1 u^2...u^{q} \in
\bigwedge\nolimits^{p} V^\ast \otimes \bigwedge\nolimits^{q} V$ has
total degree $n.$
\end{remark}

Let us define the solutions of the Maurer-Cartan equation of total
degree 1 in $B.$ The correction of   $(p+q -2)$ takes into account
that elements entering the Maurer-Cartan equation are no longer
necessarily from $B^1$ (that is when $p+q = 3$) but from $B^{p+q-2}$
for any $p$ and $q.$

\subsection{ $L_\infty$ algebra structure from the iterated adjoint action}

Let us  recall the definition of  $L_\infty$-algebras.

\begin{dfn}\label{shuffle} Consider a map
$\lambda_k \in Hom(\bigwedge\nolimits^kV, V).$  This map $\lambda_k$
acts on $\bigwedge\nolimits^nV$ for any $n \in \N$  by  coderivation
of the unshuffle coproduct on the algebra of exterior powers of $V:$
\begin{equation}\label{coder1}
 \lambda_k (v_1\wedge  \cdots \wedge  v_n) = \left\{
 \begin{array}{ll}
 0,& \ \mbox{if} \
 k > n\\
\sum_{\sigma \in Sh_n^k} (-1)^{sgn \sigma} \lambda_k(
v_{\sigma(1)}\wedge \cdots \wedge v_{\sigma(k)}) \wedge
v_{\sigma(k+1)} \cdots \wedge v_{\sigma(n)}, & \mbox{otherwise.}
 \end{array}
 \right.
 \end{equation}
The set $Sh_n^k$ is the set of all  $k$-unshuffles in  the
permutation
 group of n elements, that is all permutations such that
 $ \sigma(1) < \sigma(2)< \cdots  < \sigma(k)$ and
$ \sigma(k+1) <  \cdots  < \sigma(n).$
\end{dfn}
\begin{dfn} \label{L_infty} An $L_\infty$-structure on a graded space $V$ is a set of maps of
total degree  $1,$
$$ \lambda_k: \bigwedge\nolimits^k V \to V[2-k], \ k \in \N$$
such that the following generalized form of the Jacobi identity is
satisfied for any $n \geq 2:$
 \begin{equation} \label{Jacobi}
 \sum_{k = 1}^{n-1}[\lambda_k,\lambda_{n-k}] = 0.
\end{equation}
The bracket  $[\cdot, \cdot]$ here is the commutator of  $\lambda_k$
and $\lambda_{n-k}$ considered as operators on $\bigwedge\nolimits
V.$ These equations are higher Jacobi identities and can be
summarized in one equation:
\begin{equation}\label{hJacobi}
\lambda^2 = 0, \ \mbox{where} \ \lambda = \sum_{k \geq 1}
 \lambda_k.
\end{equation}
\end{dfn}

For a finite dimensional space $V$ there is the  following
isomorphism $\Hom(\bigwedge\nolimits^k V, \bigwedge\nolimits^l V)
\simeq \bigwedge\nolimits^k V^\ast \otimes \bigwedge\nolimits^l V. $
Any operator $\tau_{kl} \in \Hom(\bigwedge\nolimits^k V,
\bigwedge\nolimits^l V)$ is represented by $t_{kl} \in
\bigwedge\nolimits^k V^\ast \otimes \bigwedge\nolimits^l V$ and a
passage from $t_{kl}$ to $\tau_{kl}$ could be made explicit by
introducing the iterated adjoint action.

\begin{dfn}   Consider an element
$t_{kl} \in \bigwedge\nolimits^k \Vast \otimes \bigwedge\nolimits^l
V. $
 {\em The iterated adjoint action} of $t_{kl}$ on $\bigwedge\nolimits^nV$ is defined as follows.
For $n = k$ it defines a higher derived bracket as in \cite{AI,V}:
 \begin{equation}\label{had-def}
\tau_{kl}(v_1 \wedge \cdots \wedge  v_{k})= \had_{t_{kl}} (v_1
\wedge \cdots \wedge  v_{k})= [[\cdots[[ t_{kl}, v_1], v_2]\cdots],
v_{k}],\end{equation}
 thus defining a map $ \bigwedge\nolimits^k \Vast
\otimes \bigwedge\nolimits^l V \to Hom( \bigwedge\nolimits^k V,
\bigwedge\nolimits^l V): \ t_{kl} \mapsto \tau_{kl}$
  More generally, we define an action of $t_{kl}$
on any exterior power of $V,$ with values in  $\bigwedge\nolimits
\Vast \otimes \bigwedge\nolimits V$ as follows
 \begin{equation}\label{had}
\begin{array}{rl}
\had_{t_{kl}} &(v_1 \wedge \cdots \wedge  v_{n})\\
 = &\left\{
 \begin{array}{ll}
\qquad \underbrace{[[\cdots[[}_{k-n \ \mbox{\tiny times}}\!\!t_{kl},
v_1], v_2]\cdots], v_{n}],& \ \mbox{if} \
 k \geq n\\
\displaystyle{\sum_{\sigma \in Sh_n^k}} (-1)^{sgn \sigma}
\underbrace{[\cdots[}_{k \ \mbox{\tiny times}}\!t_{kl},
v_{\sigma(1)}], v_{\sigma(2)}]\cdots], v_{\sigma(k)}] \wedge
v_{\sigma{(k+1)}} \cdots \wedge v_{\sigma(n)}, & \mbox{otherwise.}
 \end{array}
 \right.\end{array}
 \end{equation}
The set of permutations  $Sh_n^k$ is as  in Definition
\ref{shuffle}. For the case $n = k$ one should keep in mind that
 $\bigwedge\nolimits^0 V^\ast = \field.$
\end{dfn}

It is shown already in \cite{Vo} that an $L_\infty$-algebra can be
obtained as a particular case of this iterated adjoint action.

 We now state  two theorems which are  $L_\infty$ analogues
 of results of Proposition \ref{Leibniz}.
 \begin{theorem} \emph{(see \cite{Vo})\ } \label{thm-Lie}
 Consider an element $L \in \La, \ L = \sum_{k=1}^\infty l_k,
 \ l_k \in \bigwedge\nolimits^k V^\ast \otimes V$
 of total degree $1,$ such that $[L,L] = 0. $  The set of maps
 $\lambda_k ={\emph \had_{l_k}} \in  \Hom( \bigwedge\nolimits^k V, V) $    form an
$L_\infty$-algebra.

In other words, $[L,L] = 0$ implies that $\lambda^2 = 0$ where
$\lambda = \had_L.$
 \end{theorem}

 The proof is a direct albeit tedious computation reducing the
condition on the maps $\lambda_k$ to $[L,L] = 0$ using the fact that
$\bigwedge\nolimits V$ is an Abelian Lie subalgebra of $\La$ and
also that $[L,v] \in \bigwedge\nolimits V^\ast$ for any $v \in
\bigwedge\nolimits V^\ast.$

We also use the following statement analogous to the one about the
Lie algebra morphism in Proposition \ref{Leibniz}.
\begin{dfn}
An $L_\infty$-morphism from an $L_\infty$-algebra $(V,\sum
\lambda_k)$ to a Lie algebra $(W, [\, , \, ])$ is a sequence of maps
$\phi_l: \bigwedge\nolimits^l V \to W$ of total  degree $0$ such
that
 we have the following equality, for all $n \geq 2,$
\[
\sum_{k+l = n } \phi_l \lambda_k = \sum_{k+l = n +1 }[\phi_k,
\phi_l].
\]
\end{dfn}

\begin{theorem}
Consider a  graded space $V$ with an $L_\infty$-structure given by
the  iterated  adjoint action of  $L = \sum l_k, \ \mbox{for} \  l_k
\in \bigwedge\nolimits^k V^\ast \otimes V, \ \mbox{and} \  [L,L] =
0.$ Consider the graded Lie algebra structure on $V^\ast \otimes V$
with the bracket $[\, , \, ]$ given by the natural pairing.

Then maps $$\emph{\had}_{l_{k+1}}: \bigwedge\nolimits^{k} V \to
V^\ast \otimes V$$ define an $L_\infty$ morphism  $ (V, \sum
\emph{\had}_{l_{k}}) \to (V^\ast \otimes V, \ [ \, ,\,  ]).$
\end{theorem}

Checking this proposition amounts just to writing the higher Jacobi
identities (\ref{Jacobi}) on $V.$
\begin{remark} In \cite[Corollary~2]{B} it is noticed that
an odd self-commuting element from $B$ defines an $L_\infty$
structure on $\bigwedge\nolimits V$ (as well as on
$\bigwedge\nolimits V^\ast$) by the iterated adjoint action.
\end{remark}

\subsection{$L_\infty$-coalgebras}
\begin{dfn}
An $L_\infty$ coalgebra structure on $V$ is a sequence of maps
$$\delta_{p}: V \to \bigwedge\nolimits^p V\ [2-p],  \ l \geq 1$$ of total
degree $1$ such that
\[ \sum_{ p \geq 1, q \geq 1} \delta_{p} \delta_{q} = 0.
\] In other words,
$\delta^2 = 0, \ \mbox{where} \
 \delta = \sum_{ p \geq1} \delta_{p}.$
\end{dfn}
Consider the Lie algebra   ${\mathcal C} = \oplus_{q \geq 0} \
V^\ast \otimes \bigwedge\nolimits^{q+1} V,$ governing the Lie
coalgebra structure on $V.$ The adjoint action of  $c_p \in V^\ast
\otimes \bigwedge\nolimits^{p} V$ on $\bigwedge\nolimits^n V$ is a
map from $\bigwedge\nolimits^n V$ to $ \bigwedge\nolimits^{n+p-1}
V:$
$$ ad_{c_p}(v_1 \wedge \cdots \wedge v_n) = \sum_{i = 1}^n
(-1)^{sgn} [c_p, v_i]\wedge v_1 \wedge \cdots \breve{v}_i \cdots
\wedge v_n.$$

\begin{theorem} \label{thm-co} Let $C = \sum_{p=1}^\infty c_p,  \ \mbox{for} \ \
c_p \in V^\ast \otimes \bigwedge\nolimits^{p} V$  be such that
$[C,C]= 0.$ Then the adjoint action of $C$ defines an
$L_\infty$-coalgebra structure on $V.$

Moreover, the iterated adjoint action of $C$ on $\bigwedge\nolimits
V^\ast$ defines an $L_\infty$-algebra structure on $V^\ast,$ the
dual space to $V.$
\end{theorem}
\begin{proof} The first statement follows from the  definition of an
$L_\infty$ coalgebra, while an $L_\infty$ structure on $V^\ast$ is a
consequence of Theorem \ref{thm-Lie}.
\end{proof}

\subsection{$L_\infty$-bialgebras}

Following the general philosophy of structures given by solutions of
Maurer-Cartan equations on differential graded Lie algebras, we
consider  solutions of an equation $[Q, Q] = 0$ where $Q$ belongs to
Lie subalgebras of $B: \ \mathcal C, \mathcal L, \mathcal B  \
\mbox{and} \ \mathcal{QB}$ from Proposition \ref{zoo}. They
respectively define an $L_\infty$-coalgebra, $L_\infty$-algebra,
$L_\infty$-bialgebra and $L_\infty$-quasibialgebra structures on
$V$. In all cases we consider the iterated adjoint action
(\ref{had-def}) of the corresponding $Q$ on $\bigwedge\nolimits V.$

Let us give two new definitions: for a differential Lie bialgebra
structure and for  a $L_\infty$ bialgebra structure.

\begin{dfn}  An
$L_\infty$ bialgebra structure on $V$ is a set of maps
$$\tau_{kl}: \bigwedge\nolimits^k V \to \bigwedge\nolimits^l V[3-(k+l)], \ k,l \geq 1$$ of total
degree $1$ such that each $\tau_{kl}=  \had_{t_{kl}}, \ t_{kl} \in
\bigwedge\nolimits^k V^\ast \otimes \bigwedge\nolimits^l V $ and
\begin{equation}\label{vse}
\sum_{k+k' =p+1} \ \sum_{l+l' =q+1} [t_{kl}, t_{k'l'}] = 0,
\end{equation}
for all $ p \geq 2, \ q \geq 2.$
\end{dfn}

In other words, an element of $\B$ of total degree $1:$
\begin{equation}\label{T}
T = \sum_{k \geq 1, \ l \geq 1} t_{kl}, \qquad t_{kl} \in
\bigwedge\nolimits^k V^\ast \otimes \bigwedge\nolimits^l V \
\mbox{such that} \ [T,T] = 0.
\end{equation} defines an
$L_\infty$-bialgebra structure on $V.$

We see that  $T $ lives  in the following bi-graded space:

\vspace{5mm}

\begin{tabular}{c|ccccc}

  %  after  \\:  \hline or \cline{col1-col2} \cline{col3-col4} ...
& $\cdots$  &     &     &    & \\
   & &     &     &    & \\
 4 \qquad & $ \bigwedge\nolimits^4 \Vast   \otimes V$     & $\cdots$     &     &  &   \\
   & &     &     &    & \\
 3 &$  \bigwedge\nolimits^3  \Vast  \otimes V$   &  $\bigwedge\nolimits^3  V^\ast \
 \ \otimes \ \bigwedge\nolimits^2 V  $  &  $\cdots$ &  &\\
  &&     &     &   &  \\
 2& $\bigwedge\nolimits^2 \Vast \otimes V$  & $\bigwedge\nolimits^2 V^\ast\ \otimes \ \bigwedge\nolimits^2 V $ &
  $\bigwedge\nolimits^2 V^\ast\ \otimes \ \bigwedge\nolimits^3 V $ &$\cdots$  &\\
   &&     &     &   &  \\
 1 & $ \Vast \otimes V$ & $  V^\ast\ \otimes \ \bigwedge\nolimits^2 V$ &
 $ V^\ast\ \otimes \ \bigwedge\nolimits^3 V$ &
 $ V^\ast\ \otimes \ \bigwedge\nolimits^4 V$ & $\cdots$ \\
   & &     &     &    & \\ \hline
 \ $k \quad \diagup \quad l $&\qquad 1 &\qquad 2  &  \qquad 3 &  \qquad 4 & \\

\end{tabular}

\vspace{5mm} Let us look at the first few equations from
(\ref{vse}).
 The first one is \begin{equation}[t_{11},t_{11}] =0,
 \label{11}
 \end{equation} providing  the equation
which defines a differential $d = \ad_{t_{11}}$ on $V.$ The
equations for $p =1, \ q = 2$ and $p =2, \ q = 1$ give respectively:
\begin{equation}\label{12 and 21}
[t_{11}, t_{12}] = 0, \qquad [t_{11}, t_{21}] = 0,
\end{equation}
 which gives the  condition that $d$ be a
  derivation of the cobracket and of the
bracket. The equation  for  $p =2, \  q = 2:$
\begin{equation}\label{22}
[t_{11}, t_{22}] +  [t_{12}, t_{21}] = 0
\end{equation}
shows that the bracket is a cocycle with respect to the cobracket up
to homotopy. The homotopy is  given by the element $t_{22} \in
\bigwedge\nolimits^2 V^\ast \otimes \bigwedge\nolimits^2 V.$

The Jacobi identity for the bracket $\had (t_{21})$ also holds only
up to homotopy; it is given by the equation for $p=3, q= 1:$
\begin{equation} \label{31}
[t_{11}, t_{31}] + \frac{1}{2} [t_{21}, t_{21}] = 0.
\end{equation}
The co-Jacobi identity up to homotopy is the equation for $p=1,q=3:$
\begin{equation} \label{13}
[t_{11}, t_{13}] + \frac{1}{2} [t_{12}, t_{12}] = 0. \end{equation}

\begin{remark} \label{L-inf bi} \begin{enumerate}
\item The homology of  $V, \  H^\ast
(V, d= \ad \ t_{11}) $ is a Lie bialgebra.
\item If $V$ is not graded, that is all $V_n = \{0\}$ for all $n
\not= 0,$ then $ V = V_0$ is an ordinary Lie bialgebra. Indeed,
since $t_{kl} $ is of total degree $1,$ on an ungraded space only
terms with $k,l$ satisfying $1-(k+l-2) = 0$ survive. The result is
that
 all $t_{kl}$ except for $ t_{12}$
and $t_{21}$ have to vanish.
The Maurer-Cartan equation then gives the axioms of a
 Lie bialgebra from Definition \ref{biLie}:
{\em $\tau_{21} = \had t_{21}$} defines a bracket
and {\em $\tau_{12} = \had t_{12}$} a cobracket.

\item If all $t_{kl} = 0 $ for all $k+l \geq 3,$ then $V$ is a
differential graded Lie bialgebra (graded Lie bialgebra with a
differential compatible with the bracket and the cobracket).
\end{enumerate}
\end{remark}

\subsection{$L_\infty$-quasi-bialgebra}

To have an $L_\infty$-quasi-bialgebra structure we need to allow
terms in $\bigwedge\nolimits^{q}V.$  Hence the Maurer-Cartan
equation of the previous subsection (\ref{T}) becomes an equation on
$S = \sum_{k \geq 0, l \geq 1} t_{kl}, \ t_{kl} \in
\bigwedge\nolimits^k V^\ast\ \otimes \ \bigwedge\nolimits^l V$ which
differs from $T$ because it contains elements $t_{0l} \in
\bigwedge\nolimits^l V.$ This $S$ must satisfy a set of equations
indexed by $p \geq 1, q \geq 2:$
\[
[S,S] = \sum_{k+k' = p} \ \sum_{l+l' =q } [t_{kl}, t_{k'l'}] = 0,
\]
 \vspace{5mm}

\begin{tabular}{c|ccccc}

  %  after  \\:  \hline or \cline{col1-col2} \cline{col3-col4} ...
& $\cdots$  &     &     &    & \\
   & &     &     &    & \\
 4 \qquad & $ \bigwedge\nolimits^4 \Vast   \otimes V$     & $\cdots$     &     &  &   \\
   & &     &     &    & \\
 3 &$  \bigwedge\nolimits^3  \Vast  \otimes V$   &  $\bigwedge\nolimits^3  V^\ast \
 \ \otimes \ \bigwedge\nolimits^2 V  $  &  $\cdots$ &  &\\
  &&     &     &   &  \\
 2& $\bigwedge\nolimits^2 \Vast \otimes V$  & $\bigwedge\nolimits^2 V^\ast\ \otimes \ \bigwedge\nolimits^2 V $ &
  $\bigwedge\nolimits^2 V^\ast\ \otimes \ \bigwedge\nolimits^3 V $ &$\cdots$  &\\
   &&     &     &   &  \\
 1 & $ \Vast \otimes V$ & $  V^\ast\ \otimes \ \bigwedge\nolimits^2 V$ &
 $ V^\ast\ \otimes \ \bigwedge\nolimits^3 V$ &
 $ V^\ast\ \otimes \ \bigwedge\nolimits^4 V$ & $\cdots$ \\
   & &     &     &    & \\
   0 &
    &  &
 $\bigwedge\nolimits^3 V$ &
 $ \bigwedge\nolimits^4 V$ & $\cdots$ \\
   & &     &     &    & \\ \hline
 \ $(p+1) \quad \diagup \quad (q +1)$&\qquad 1 &\qquad 2  &  \qquad 3 &  \qquad 4 & \\

\end{tabular}
%\hbox to 8pt {$V$\hss xx}\hbox to 18pt {$\bigwedge\nolimits^2 V$\hss {xxx}}
\vspace{5mm}

Terms  $ t_{0l}$ do not act on $\bigwedge\nolimits V,$ however they
change the co-Jacobi condition and the other equations as well in
comparison with the equations on $T$.

In particular, terms $t_{01} \in V$ and $t_{02} \in V \wedge V$
change the nature of certain  equations: for example, $[t_{11},
t_{11}] = 0$ is no longer true in the presence of $t_{01}.$ Allowing
non-zero terms $t_{01}$ and $t_{02}$  would lead us to a completely
different setup of weak $L_\infty$-algebras and to avoid that we
impose that these terms are $0.$ Weak $L_\infty$-algebras appear in
physics, vanishing of the terms $t_{01}$ and $t_{02}$ is related to
certain boundary condition in Batalin-Vilkovisky formulation of the
open string A-model (see \cite{OP}).

\begin{dfn}Consider an element of $\QB$ of total degree $1:$
\begin{equation}\label{0}
S = \sum_{k \geq 1, l \geq 1} t_{kl} + \sum_{l\geq3} t_{0l}, \
t_{kl} \in \bigwedge\nolimits^k V^\ast\ \otimes \
\bigwedge\nolimits^l V,
\end{equation} such that
\begin{equation}\label{bracket S}
[S,S] = 0.
\end{equation}
An $L_\infty$ quasi-bialgebra structure on $V$ is the set of maps
$$\tau_{kl}: \bigwedge\nolimits^k V \to \left(\bigwedge\nolimits^l V\right)[3-(k+l)], \ k \geq 0, \ l \geq 1$$
of total degree $1,$ where  $ \tau_{kl} = \had t_{kl}.$
\end{dfn}

 Let us state properties
of an $L_\infty$-quasi-bialgebra similar to the properties stated in
Remark \ref{L-inf bi}:

\begin{remark} \label{L-inf_quasi-bi} \begin{enumerate}
\item The homology of $H^\ast V = H^\ast
(V, \ad t_{11}) $ is a Lie quasi-bialgebra.
\item If $V$ is not graded then $ V $ is an ordinary Lie quasibialgebra. Indeed,
since $t_{kl} $ is of degree $2-(k+l-1)$ all $t_{kl}$ except for $
t_{12},  t_{21}$ and $t_{03}$  have to vanish. Then  {\em $\had
t_{12}$} defines a cobracket, {\em  $\had t_{21}$} a bracket and
together with  $t_{03} \in \bigwedge\nolimits^3 V$ they satisfy the
axioms of Lie quasibialgebras from Definition \ref{qbiLie}.
\item If all $t_{kl} = 0 $ for all $k+l \geq 3$ then $V$ is a
differential graded Lie quasi-bialgebra.
\end{enumerate}
\end{remark}

\begin{remark} \label{thm} For $L_\infty$-(quasi)bialgebras
there is no analogue of Theorems \ref{thm-Lie}, \ref{thm-co}. For $Q
\in B$ The condition $[Q,Q] = 0$ in general does not imply  {\em $
(\had Q)^2 = 0$}.   Although it is true for $Q$ either in $\mathcal
L$ or in $\mathcal C,$ in a more general case the equation {\em $
(\had Q)^2 = 0$} force higher order conditions on $Q$  which are not
necessarily satisfied  even for an ordinary (not $L_\infty$) Lie
bialgebra or a Lie quasi-bialgebra. For example for a Lie bialgebra
the condition $ (\had Q)^2 = 0$ would force an extra axiom. Namely,
the cobracket applied to an element and then the bracket to the
result must be required to be $0,$ which gives a "no cycle"
condition but it is not in the initial set of axioms for a Lie
bialgebra. That kind of bialgebras "without cycles" appears in the
work of Chas and Sullivan on string topology \cite{CS}.
 \end{remark}

\section{Formal geometry of $L_\infty$-structures. }

%%%%%%%%%%%%%%%%%%

We now move to  a geometric definition of an $L_\infty$-algebra
structure on a finite dimensional space (see for example \cite{AKSZ}
or \cite{M}). Consider a finite dimensional graded vector space $V=
\oplus_k V_k$ with a chosen basis $\{e_\alpha\}.$ Then a Lie algebra
structure is defined by structure constants $c_{\alpha,
\beta}^\gamma:$
\[
[e_\alpha, e_\beta] = \sum c_{\alpha, \beta}^\gamma e_\gamma,
\]
satisfying the structure equation which boils down to the Jacobi
identity. We could also see the bracket acting from  $V \wedge V$ to
$V$  as a differential operator $Q_{21} = c_{\alpha,\beta}^\gamma
e_\gamma \otimes \displaystyle{\frac{\partial}{\partial e_\alpha}
\wedge \frac{\partial}{\partial e_\beta}}.$  We use the index $21$
for this operator to underline
that it is  quadratic in $\displaystyle{\frac{\partial}{\partial
e}}$'s and linear in $e$'s.

Under the identification similar to (\ref{sym})
\begin{equation}\label{sym2} \ \Sym^n (X[1]) \simeq
\left(\bigwedge\nolimits^n (X)\right)[n]\end{equation} a Lie bracket
$\bigwedge\nolimits^2 V \to V$ is a map of degree 1: if we shift the
degree on the right and on the left by $2,$
$\left(\bigwedge\nolimits^2V\right)[2] \to V[2],$ we get
equivalently $Sym^2 (V[1]) \to (Sym (V[1]))[1].$ On the other hand,
the dual of the free cocommutative coalgebra $Sym(V[1])$ could be
identified with the algebra of formal power series  on the
$\Z$-graded space $V[1].$  We say that the algebra of formal power
series  defines a formal manifold. This way a Lie bracket becomes a
degree $1$ quadratic vector field on it. The dual space to $V[1]$ is
$V^\ast[-1].$ Elements of $Sym(V^\ast[-1])$ are coordinate functions
on $V[1].$ Let us chose corresponding coordinates $\{x^\alpha \},$
dual to $e_\alpha,$ such that  in the space $V^\ast[-1] \ x^\alpha$
has degree shifted by $1: \ \tilde{x}^{\alpha} + 1.$

A very natural generalization of this construction is the definition
of an $L_\infty$-algebra structure:

\begin{proposition} (see for example \cite{K})
An $L_\infty$-algebra structure on a graded space $V$ is given by
 vector field $Q$ of degree $1$ on the $\Z$-graded formal
manifold corresponding to $V[1]$ such that $[Q,Q] = 0.$
\end{proposition}
\begin{dfn}(following Vaintrob \cite{V})
 A non-zero self-commuting operator (and in particular a
vector field) on a formal manifold is called homological.
\end{dfn}
In particular an $L_\infty$-algebra structure is given by a
homological vector field:
\begin{eqnarray*}
 Q =  \displaystyle{b_\alpha^\beta x^\alpha
\frac{\partial}{\partial x^\beta}+ c_{\alpha \beta}^\gamma x^\alpha
x^\beta\frac{\partial}{\partial x^\gamma} + d_{\alpha \beta
\gamma}^\delta x^\alpha x^\beta x^\gamma \frac{\partial}{\partial
x^\delta} + \cdots, }\\
Q = \displaystyle{\sum_{k=1}^\infty Q_{k1}  , \qquad Q_{k1} =
f_{\alpha_1 \cdots \alpha_k}^\beta x^{\alpha_1} \cdots
x^{\alpha_k}\frac{\partial}{\partial x^\beta}}.
\end{eqnarray*}
 In fact, each term in $Q$ corresponds to a map
of degree $1$ on the dual space:
\begin{equation} \label{Q}
Q = \sum Q_{k1}, \  Q_{k1}: V^\ast[-1] \to \left( Sym^k
(V^\ast[-1])\right)[1].
\end{equation}

Counting degrees after the identification (\ref{sym2}) we get
$Q_{k1}: V^\ast[-1] \to \left(\bigwedge\nolimits^kV^\ast[-k]\right)
[1]$ which on the dual space becomes a map $\bigwedge\nolimits^kV
\to V[2-k]$ as in Definition \ref{L_infty}.

This way we are naturally brought  to a geometric definition of an
$L_\infty$-bialgebra  structure on $V.$ In the sum (\ref{Q}) we
should consider not only $Q_{k1}: V^\ast \to Sym^k V^\ast,$ but any
$Q_{km}: Sym^m V^\ast \to Sym^k V^\ast.$ It becomes  a differential
operator on $V[1],$ acting  on functions on  $V[1],$ which are
elements of  $Sym (V^\ast[-1]).$

According to Definition \ref{degree} the operator $Q_{km}$ has the
external degree $\textbf{ed}(Q_{km}) = k+m -2,$ following the
external grading in the space $ \Lambda^m V^\ast\otimes \Lambda^k
V.$ The total degree of $Q_{km}$ is then a difference of internal
degrees of the result with the initial element  plus the external
degree $k+m-2.$

We consider the ring of differential operators on $V[1]:$
$$\A(V): = Sym(V^\ast[-1]) [\xi^1, \cdots, \xi^n]$$ where
$\xi^\alpha =  \displaystyle{\frac{\partial}{\partial x_\alpha}}.$
Notice that  $ \xi^\alpha \in \left(V[1]\right).$ On $\A(V)$ there
is a standard graded Poisson bracket of \cite{BV,FV}. For any two
elements $E, F \in \A(V)$ it could be defined in coordinates as
follows:
\[
\{E,F\} =  \frac{\partial E}{\partial x^\alpha} \frac{\partial
F}{\partial \xi_\alpha}  - (-1)^{\tilde{E}\tilde{F}}\frac{\partial
F}{\partial x^\alpha}
 \frac{\partial E}{\partial \xi_\alpha}
\]

\begin{proposition}
Any $L_\infty$-bialgebra structure on $V$ is given by an operator of
degree $1$ on the corresponding formal manifold. The Poisson bracket
of this operator with itself is $0.$
\end{proposition}
\begin{proof}
Consider an operator $Q = \sum Q_{pq}, \ \mbox{where} \
 Q_{pq} = v^1\cdots  v^{p} w^\ast_{1}\cdots w^\ast_{q}.$
Being of degree $1$ gives the following equality: $\sum
(\widetilde{v}_i + \widetilde{w}_j ) +p+q - 2= 1.$ Here
$\widetilde{w}_j = - \widetilde{w^\ast}_j.$ The operator $Q$ acts on
skew-symmetric algebraic functions on $V,$ that is on
 $Sym(V^\ast[-1]).$ Thus, after the
identification (\ref{sym}) we get that $Q_{pq}: \bigwedge\nolimits^p
V^\ast \to \left(\bigwedge\nolimits^q V^\ast\right)[1][2 -(p+q)]$
which on the dual space gives $ \bigwedge\nolimits^q V [p+q -3] \to
 \bigwedge\nolimits^p V$ leading to the right degrees.

The Poisson bracket of $Q$ with itself  being $0$ corresponds to the
self-commuting condition in (\ref{T}).
\end{proof}

\section{Manin $L_\infty$-triple and $L_\infty$-pair}

We give a definition of a Manin $L_\infty$-triple (which also could
be called a strongly homotopy Manin triple), a natural
generalization of a Manin triple.
 For this we need a notion of an $L_\infty$-subalgebra:
\begin{dfn}
Let $(W, Q)$ be an $L_\infty$-algebra. Then a subspace $V \subset W$
is an $L_\infty$-subalgebra if $V$ is $Q$-invariant.
\end{dfn}
It means that $V$ is an $L_\infty$-subalgebra of $W$ if the image of
the operator $Q$ acting on $V^\ast[-1]$ is in $Sym V^\ast[-1][1].$
In other words it means that  $Q= \sum_{k,l \geq 1} Q_{k1}$
restricted to $V$ consists of maps $ V ^\ast \to
\bigwedge\nolimits^k V^\ast[2-k]$ or, on the dual, it gives an
$L_\infty$ structure on $V, \ Q_{k1}: \bigwedge\nolimits^k V \to V
[2-k].$

\begin{dfn} A finite dimensional Manin  $L_\infty$-triple is a triple
of finite dimensional $L_\infty$-algebras $(\g, \g_+, \g_{-})$
equipped with a nondegenerate bilinear form $<,>$ such
that
\begin{itemize}
    \item $\g_+, \g_{-}$ are $L_\infty$-subalgebras of $\g$ such
    that $\g =  \g_+ \oplus \g_{-}$ as a vector space;
    \item $\g_+$ and $\g_{-}$ are isotropic with respect to $<,>;$
    \item the $n$-brackets constituting the $L_\infty$-algebra structure
$\lambda_n: \bigwedge\nolimits^n (V^\ast \oplus V) \to V^\ast \oplus
V$ are invariant with respect to the bilinear form $<,>,$ that is
\begin{equation} \label{invariance}
<\lambda_n(v_1 \wedge \cdots \wedge  v_n), v_0> =
(-1)^{\widetilde{v_n} \widetilde{v_0}}
<\lambda_n(v_1 \wedge \cdots \wedge v_{n-1} \wedge v_0), v_n>
\end{equation}
\end{itemize}
\end{dfn}
Notice that a Manin triple is an example of a Manin
$L_\infty$-triple.

The invariance is in fact cyclic, since $\lambda_n$ (being a map on
$\bigwedge\nolimits^n V$) is antisymmetric in all variables.

\begin{theorem}\label{Manin triple}
The notions of a (finite dimensional)  Manin $L_\infty$-triple and
an $L_\infty$-bialgebra are equivalent.
\end{theorem}
\begin{proof}
Consider a Manin $L_\infty$-triple $(\g, \g_+, \g_{-}).$ Let us show
that there is an $L_\infty$-bialgebra structure on $\g_+.$ Let the
$L_\infty$-algebra structure on $\g$ be given by maps $\lambda_k:
\bigwedge\nolimits^k \g \to \g,$ that is  $\lambda_k:
\bigwedge\nolimits^k (\g_+ \oplus \g_-) \to (\g_+ \oplus \g_-).$
These maps can be split as follows:
$$\lambda_{mn+}: \bigwedge\nolimits^m \g_+ \otimes \bigwedge\nolimits^n \g_- \to \g_+ \ \mbox{and} \
 \lambda_{mn-}:\bigwedge\nolimits^m \g_+ \otimes \bigwedge\nolimits^n \g_- \to \g_-
\ \mbox{ with} \ m+n = k.$$

Using non-degeneracy and invariance  (\ref{invariance}) of the
internal product $<,>$ after an identification of $\g_-^\ast$ with
$\g_+,$ we get the equivalences
$$\begin{array}{l}
 \Hom_{inv} (\bigwedge\nolimits^m \g_+ \otimes
\bigwedge\nolimits^n \g_-,\g_+) \simeq \Hom(\bigwedge\nolimits^m
\g_+, \bigwedge\nolimits^{n+1} \g_+)\simeq \bigwedge\nolimits^{m}
\g_+^\ast \otimes \bigwedge\nolimits^{n+1} \g_+ \quad \mbox{and} \\
  \Hom_{inv} (\bigwedge\nolimits^m \g_+ \otimes \bigwedge\nolimits^n \g_-,\g_-)\simeq
\Hom (\bigwedge\nolimits^{m+1} \g_+, \bigwedge\nolimits^{n} \g_+)
\simeq \bigwedge\nolimits^{m+1} \g_+^\ast \otimes
\bigwedge\nolimits^{n} \g_+.\end{array}$$

Let  $l_{mn+} \in \bigwedge\nolimits^{m} \g_+^\ast \otimes
\bigwedge\nolimits^{n+1} \g_+$ correspond to $\lambda_{mn+}.$ This
way, $l_{mn+}$ would act by a iterated adjoint action:
\[
\had(l_{mn+}): \bigwedge\nolimits^m \g_+ \to
\bigwedge\nolimits^{n+1} \g_+.
\]

On the other hand $\lambda_{mn-}:\bigwedge\nolimits^m \g_+ \otimes
\bigwedge\nolimits^n \g_- \to \g_-$ again is given by a iterated
adjoint action of an element  $l_{mn-} \in \bigwedge\nolimits^{m+1}
\g_+^\ast \otimes \bigwedge\nolimits^{n} \g_+$  and
 \[
\had(l_{mn-}): \bigwedge\nolimits^{m+1} \g_+ \to
\bigwedge\nolimits^{n} \g_+.
\]
Now, $L = \sum_{m,n \geq 1} (l_{mn+} + l_{mn-})$ satisfies $[L,L]=
0$ from the condition on $\lambda.$

In the other direction it is actually easier.
 Let an $L_\infty$-bialgebra structure on $V$ be defined by the
iterated adjoint action of $$T = \sum_{k \geq 1, l \geq 1} t_{kl}, \
\mbox{ where} \ t_{kl} \in \bigwedge\nolimits^k V^\ast \otimes
\bigwedge\nolimits^l V.$$ Then the corresponding Manin
$L_\infty$-triple is $(V \oplus V^\ast, V, V^\ast)$ with the
nondegenerate form given by the natural pairing of $V$ and $V^\ast.$

The $L_\infty$ structure on $V^\ast \oplus V$ is given by $\tau =
\had(T).$ In
$$\had(t_{kl}): (V^\ast \oplus V) \to \bigwedge\nolimits^{k-1} V^\ast \otimes
\bigwedge\nolimits^l V \oplus \bigwedge\nolimits^k V^\ast \otimes
\bigwedge\nolimits^{l-1} V$$ we recognize the $L_\infty$-coalgebra
structure dual to an $L_\infty$-algebra structure we are looking
for.

To get the $L_\infty$-subalgebra structure on $V$ we consider the
quadratic equations (\ref{vse}) for $q = 2.$ We get:
\[
\sum_{k+k'= p} t_{k1} t_{k'1} = 0,
\]
the set of  equations defining  an $L_\infty$-structure
 on $V$ (by taking $ t_{k1} = l_k$ in (\ref{thm-Lie})).

The operator  $\tau_V = \sum_{k \geq 1} \had(t_{k1}): V^\ast \to
\bigwedge\nolimits V^\ast$ is the restriction of $\tau$ to $V^\ast$
making $V$ an $L_\infty$-subalgebra of $V^\ast \oplus V,$ since
$T_V$ is self-commuting.

In the same way, elements $ t_{1l} \in V^\ast \otimes
\bigwedge\nolimits^lV$ satisfy the following quadratic equation
(equations (\ref{vse}) for $p = 2$):
\[ \sum_{l+l'= q} t_{1l} t_{1l'} = 0
\]
therefore they induce maps $\had(t_{1l}): V \to \bigwedge\nolimits^l
V,$ on $V$ thus defining  an $L_\infty$-algebra structure on
$V^\ast.$
\end{proof}

 The notion which
leads to an $L_\infty$-quasi-bialgebra is a Manin $L_\infty$-pair
(strongly homotopy Manin pair):
\begin{dfn} A finite dimensional Manin $L_\infty$-pair is a pair
of finite dimensional $L_\infty$-algebras $\g \supset \g_+$ equipped
with a nondegenerate invariant (\ref{invariance}) bilinear form
$<,>,$ such that
 $\g_+$ is an isotropic $L_\infty$-subalgebra of $\g.$

 A Manin $L_\infty$-quasi-triple is a Manin $L_\infty$-pair $(\g,\g_+)$ with a
 chosen Lagrangian complement of $\g_+.$
\end{dfn}
There is also a correspondence as in Theorem \ref{Manin triple} with
a similar proof:
\begin{theorem}
The notions of a (finite dimensional) Manin $L_\infty$-quasi-triple
and an $L_\infty$-quasi-bialge\-bra are equivalent.
\end{theorem}

\end{document}